\documentclass[12pt]{article} 
\usepackage[margin=0.7in]{geometry}                           

\usepackage{authblk}
\usepackage{bm}
\usepackage{amsmath}
\usepackage{amsthm}
\usepackage{amssymb}
\usepackage{graphicx}
\usepackage[font=footnotesize]{caption}
\usepackage{float}

\newtheorem{thm}{Theorem}[section]

\newtheorem{defn}[thm]{Definition}

\newtheorem{prop}[thm]{Proposition}

\date{}                                                        

\begin{document}

\title{Topological Signals of Singularities in Ricci Flow}

\author[1]{\large Paul M. Alsing\footnote{paul.alsing@us.af.mil}}
\author[2]{Howard A. Blair\thanks{blair@ecs.syr.edu}}
\author[3]{Matthew Corne\thanks{cornem@uwstout.edu}}
\author[4]{Gordon Jones\thanks{gojones@syracuse.edu}}
\author[5]{\\Warner A. Miller\thanks{wam@fau.edu}}
\author[6]{Konstantin Mischaikow\thanks{mischaik@math.rutgers.edu}}
\author[7]{Vidit Nanda\thanks{vnanda@sas.upenn.edu}}
\affil[1]{\normalsize Air Force Research Laboratory, Information Directorate, Rome, NY 13441}
\affil[2]{Department of Electrical Engineering and Computer Science, Syracuse University, Syracuse, NY 13244}
\affil[3]{Department of Mathematics, Statistics, and Computer Science, University of Wisconsin--Stout, Menomonie, WI 54751}
\affil[4]{Department of Electrical Engineering and Computer Science, Syracuse University, Syracuse, NY 13244}
\affil[5]{Department of Physics, Florida Atlantic University, Boca Raton, FL 33431}
\affil[6]{Department of Mathematics, Rutgers University, Piscataway, NJ 08854}
\affil[7]{Department of Mathematics, University of Pennsylvania, Philadelphia, PA 19104}

\maketitle

\begin{abstract}

We implement methods from computational homology to obtain a topological signal of singularity formation in a selection of geometries evolved numerically by Ricci flow. Our approach, based on persistent homology, produces precise, quantitative measures describing the behavior of an entire collection of data across a discrete sample of times. We analyze the topological signals of geometric criticality obtained numerically from the application of persistent homology to models manifesting singularities under Ricci flow. The results we obtain for these numerical models suggest that the topological signals distinguish global singularity formation (collapse to a round point) from local singularity formation (neckpinch). Finally, we discuss the interpretation and implication of these results and future applications.

\end{abstract}



\section{Introduction}
\label{intro}
Ricci flow (RF) \cite{Hamilton:1982} is a system of partial differential equations (PDEs) that has been used in the classification of two- and three-dimensional geometries \cite{CK,Thurston:1997} via the uniformization theorem and geometrization conjecture, respectively. An interesting application involves the classification of singularities in the solutions of the RF equations. In two dimensions, compact surfaces with two-sphere topology collapse to round points, while in three dimensions, a finite number of pinching singularities can occur in different parts of the geometry, depending on the problem being considered. Such singularity formation has been thoroughly studied (cf. \cite{AIK:2011,AIK:2015,AK,GZ}). Outside of purely theoretical considerations, it is useful in physics, notably in its relation to the renormalization-group (RG) flow in string theory 
\cite{Perelman:2003-1,Perelman:2003-2,Carfora:2010,Carfora:2011,Carfora:2017,vanRaamsdonk,Woolgar}, and in discrete analogues used in imaging 
\cite{Gu:2012-1,Gu:2012-2,Miller:2014}.

In view of the applied problems, solutions to the RF equations are represented as datasets obtained via numerical integration. Now, suppose that one is given a collection of sparse datasets, without a~metric, representing an object altered under RF. This poses a challenge if one desires to know a metric or the curvature of the object at different times. While different interpolative techniques exist, they will not yield a unique answer to a given problem. Nonetheless, it is worthwhile to explore whether~a~common feature is shared in the evolution of such datasets because, if so, it would allow one to identify the type of system without having to know information about the metric. In particular, collapse of a topological two-sphere to a round point represents a global singularity, while the neckpinching singularities of initially-pinched topological three-spheres are local singularities. Given such datasets, then, one can seek a~topological characterization to indicate whether the datasets represent global or local singularities.

Persistent homology (PH) \cite{C,EH:2009,KMM,NS} is a method of topological data analysis (TDA) useful for obtaining such a characterization. One uses PH by assigning data to appropriate complexes (e.g.,~simplicial, \u{C}ech, Vietoris--Rips). Then, information about the Betti number spectra is used to characterize systems by filtering according to threshold values and assessing emergent features of the collection of data. Applications of PH include the analysis of force networks in granular media \cite{KGKM:2014}, robotic path planning \cite{BGK}, sensor network coverage \cite{dSG}, target tracking \cite{Bendich} and cosmic structure~\cite{vdWeygaert}. Persistent homology provides information about the global structure of the data. However, it also provides a detailed representation of the local features of the dataset by carrying geometric information about each of the elements.

Numerical simulations typically limit us to pointwise snapshots at single moments of time of quantities of interest, even with a complete set of data and methods for collection and measurement. Further, such simulations are plagued with stiffness that obstructs or distorts the evolution of a system. Persistent homology provides a method for evaluating a collection of and tracking trends in data across an entire domain of interest and at each step of the time evolution by identifying whether multiple small effects or fewer large effects dominate changes in the geometry. 

In this work, we apply persistent homology to numerical models of Ricci flow of a selection of metrics on $ S^{2} $ and $ S^{3} $. Persistent homology captures a topological signal of the formation of critical geometric phenomena obtained through numerical solutions of models in RF. Specifically, this~signal consists of the Betti number spectra, tables measuring lengths of persistence intervals and~distances between persistence diagrams. This provides a quantitative method for identifying common trends in the resulting collection of data obtained by solution of these PDEs. To produce numerical simulations of RF on $ S^{2} $, we consider a surface produced by a function from a random number generator. On~$ S^{3} $, we have developed a parametric ``interpolative'' dumbbell such that one can adjust a~parameter and numerically interpolate between a symmetric dumbbell to an asymmetric dumbbell and, with~appropriate choices for the remaining parameters, a degenerate dumbbell. The~degenerate dumbbell might allow for the study of Type-II singularity formation as described in \cite{GZ}, though we do not explore this due to difficulties in obtaining such a singularity \cite{GI:2005,GI:2009}. Also on $ S^{3} $, we consider a dimpled dumbbell (multiple necks and bumps), also called the ``genetically-modified peanut'' \cite{AK}. While other authors have also considered numerical models devoted explicitly to RF on $ S^{3} $ (e.g., \cite{GI:2005,GI:2009}) and the models considered here are based on well-established theoretical results \cite{AIK:2011,AIK:2015,CK}, our construction regarding the interpolative dumbbell has not been reproduced in the literature, as far as we can ascertain. We emphasize that this work concerns the application of persistent homology to datasets representing numerical solutions of Ricci flow. We hope to use the knowledge gained in this work to consider less robust (e.g., without metric) datasets in this topological setting or to have the topological and geometric perspectives inform each other. These numerical models suggest that PH provides distinct topological signals for the formation of global singularities and local singularities. This is due to differences in the behavior of the datasets describing these models: we find that global singularities correspond with many small changes in the data, while local singularities are punctuated by a single large change in the data. While this is not unexpected, the precise signal has not been observed thus far, and it provides an opportunity for future explorations at the interface of geometry and topology, notably information science, quantum tomography and problems in network congestion.

In Section \ref{sec:Preliminaries}, we review the necessary mathematical preliminaries. These include brief expositions of persistent homology, Ricci flow and singularity formation in RF. In Section \ref{sec:Methodology}, we describe the models, our data-generating algorithm, the initial and boundary conditions for the different models, the assignment of the data to a triangulation and the algorithm used to acquire the topological signals of the models. In Section \ref{sec:Results}, we present the analysis of the models. Our results show that PH provides distinct topological signals for the models when uniformization occurs, when collapse occurs or when a singularity is encountered. We conclude with Section \ref{sec:Conclusions}, where we review the results and discuss their interpretation and implication. Finally, we present future applications.

\section{Preliminaries and Background}
\label{sec:Preliminaries}

We begin with summaries of {persistent homology} and {Ricci flow}. Since our treatment of these topics is cursory, we refer the interested reader to \cite{C,EH:2009,Ghrist,NS} for overviews of persistent homology and its applications and to \cite{Miller:2014crf,Miller:2014} for various technical details pertaining to Ricci flow and its discrete~analogues.

\subsection{Persistent Homology}
\label{ph}

A {filtration} consists of a one-parameter family of triangulable topological spaces $X_p$ indexed by some real number $p \in \mathbb{R}$ subject to the constraint that $X_p$ is a subspace of $X_q$ whenever $p \leq q$. At its core, persistent homology associates to each filtration a sequence of intervals $[b_*,d_*)$, called {persistence intervals} or {lifespans}, which are indexed by homology classes across all of the $X_p$s as $p$ ranges from $-\infty$ to $\infty$. Here, we describe how those intervals are obtained in the special case where each intermediate topological space in the filtration is a finite simplicial complex $\mathcal{K}_p$. For the purposes of this section, we~do not require an ambient Euclidean space containing all of the simplicial complexes; it suffices to assume finiteness and the subcomplex relation $\mathcal{K}_p \hookrightarrow \mathcal{K}_q$ whenever $p \leq q$.

\begin{defn}
\label{def:phgroup}
Given a simplicial filtration $\{\mathcal{K}_p \mid p \in \mathbb{R}\}$, the $d$-dimensional $\ell$-{{persistent homology group}} of $\mathcal{K}_p$ is the image of the map $H_d(\mathcal{K}_p) \to H_d(\mathcal{K}_{p+\ell})$ induced on (ordinary) homology groups by the inclusion of $\mathcal{K}_p$ into $\mathcal{K}_{p+\ell}$.
\end{defn}
In general, even the space of simplicial filtrations is too large for persistent homology groups to provide a useful invariant, so one adds a local finiteness assumption as follows. A simplicial filtration $\{\mathcal{K}_p \mid p \in \mathbb{R}\}$ is said to be {tame} whenever the set of indices $p$ for which the inclusion $\mathcal{K}_{p-\epsilon} \hookrightarrow \mathcal{K}_{p+\epsilon}$ does {not} induce an isomorphism on homology for all dimensions is finite. Henceforth, we will assume that all filtrations are tame.

One may compute homology with coefficients in a field to obtain, for each simplicial filtration $\{\mathcal{K}_p\}$ and for each dimension $d \geq 0$, a one-parameter family of vector spaces $\{H_d(\mathcal{K}_p)\}$, which comes equipped with linear maps $\phi_{p\to q}: H_d(\mathcal{K}_p) \to H_d(\mathcal{K}_q)$ whenever $p \leq q$. These maps are induced on homology by including $\mathcal{K}_p$ into $\mathcal{K}_q$ and therefore satisfy the associative property: the composite $H_d(\mathcal{K}_p) \to H_d(\mathcal{K}_q) \to H_d(\mathcal{K}_r)$ is the same as $H_d(\mathcal{K}_p) \to H_d(\mathcal{K}_r)$ across any increasing triple $p \leq q \leq r$ of positive real numbers. Such a family of vector spaces and linear maps is called a {persistence module}.

Two key results involving persistence modules make persistent homology an ideal candidate for investigating the topological changes in a filtered space. The first has to do with a finite representation. Associated with any interval $[a,b]$ with $b \geq a$ is the {interval module} $\mathbf{I}_{[a,b]}$ whose constituent vector spaces are trivial whenever $p \notin [a,b]$ and one-dimensional otherwise; the linear maps are identities whenever the source and target vector spaces are both nontrivial (we allow semi-infinite intervals here, so the right endpoint is allowed to equal $\infty$). The following proposition was first proved in \cite{ZC} and has been substantially generalized since \cite{CdSGO}.
\begin{prop}
\label{C2}
Every persistence module arising from a tame simplicial filtration is canonically isomorphic to a direct sum of {{interval modules}}.
\end{prop}
In particular, the isomorphism classes of persistence modules are conveniently indexed by a finite family of intervals of the form $[a,b]$. The $d$-dimensional {persistence diagram} (PD) of a tame simplicial filtration $\{\mathcal{K}_p\}$ is defined to be the collection of intervals that arise from the decomposition of its $d$-dimensional persistence module.

The second key result confers robustness to persistent homology and is usually called the {stability theorem} \cite{CdSGO}. The basic idea is that the space of all simplicial filtrations and the space of persistence diagrams (PDs) may both be assigned natural metric structures so that the process of associating to each filtration its PDs is one-Lipschitz.

\subsection{Ricci Flow}
\label{rf}

Let $ \mathcal{M} $ be a smooth compact Riemannian manifold with metric $ \mathbf{g} $ and $ \mathcal{V} (\mathcal{M}) $ the set of vector fields on $ \mathcal{M} $. Then, the {Riemann curvature tensor} is the tensor $ \mathbf{Rm}: \mathcal{V} (\mathcal{M}) \times \mathcal{V} (\mathcal{M}) \times \mathcal{V} (\mathcal{M}) \rightarrow \mathcal{V} (\mathcal{M}) $ given by:
\begin{equation}
\label{eq:Riem}
R(X,Y,Z) = \nabla_{X}\nabla_{Y}Z - \nabla_{Y}\nabla_{X}Z - \nabla_{[X,Y]}Z 
\end{equation}

\noindent where $ \nabla_{U}V $ is the covariant derivative of $ V $ with respect to the vector field $ U $ \cite{BG,Nakahara}. The Riemann tensor measures the difference between the resulting vectors obtained at a point $ q \in \mathcal{M} $ by parallel transporting a vector along two different paths starting at a point $ p \in \mathcal{M} $.

The Ricci flow (RF) is a geometric flow that evolves the metric $ \mathbf{g} $ on $ \mathcal{M} $ by the equation:
\begin{eqnarray}
\label{eq:RF}
\frac{d}{dt}\mathbf{g} = -2\,\mathbf{Ric},
\end{eqnarray}

\noindent where the differentiation on the left--hand side is with respect to an external time parameter, $ t $. The {Ricci tensor} $ \mathbf{Ric} $ is a contraction of the Riemann tensor $ \mathbf{Rm} $, producing a tensor of type $ (0,2) $ from a tensor of type $ (1,3) $, 
and provides a certain average of the sectional curvatures of all two-planes 
along a given direction. 
Ricci flow was introduced in the early 1980s by Richard S. Hamilton \cite{Cao:2003,Hamilton:1982}.

The RF Equation (\ref{eq:RF}) yields a forced diffusion equation for the scalar curvature $ R $ so that it evolves~as:
\begin{equation}
\label{eq:fdeqn}
\dot R = \triangle R + 2\, |\mathbf{Ric}|^2
\end{equation}

\noindent and tends to uniformly distribute curvature over the manifold. Here, $\triangle_{\mathbf{g}}$ is the usual Beltrami--Laplace operator with respect to the metric, and $ | \cdot | $ is the tensor norm. 
 
In three and higher dimensions, RF can develop singularities in curvature. For three-dimensional manifolds, such flows can admit a finite number of these singularities at distinct times in the interval $t \in [0,\infty)$. Identifying, understanding and surgically removing these singularities was an essential step in Hamilton's program and Perelman's proof of Thurston's geometrization conjecture to decompose a~manifold into its prime and simple pieces \cite{CK}.


We review the work on singularities in RF found in \cite{CK,Hamilton:1982,Hamilton:1995}. Consider a solution $ (\mathcal{M},\mathbf{g}(t)) $ of RF that exists on a maximal time interval $ [0, T) $. Such a solution is {maximal} if $ |\mathbf{Rm}| $, where $ \mathbf{Rm} $ is the Riemann tensor, is unbounded as $ t \rightarrow T $, $ T < \infty $ or $ T = \infty $. If $ T < \infty $ and $ |\mathbf{Rm}| $ becomes unbounded as $ t \rightarrow T $, then a maximal solution develops singularities, and $ T $ is called the {singularity time} \cite{GZ}.

There exist several types of singularities.

\begin{defn}
\label{def:singularities}
\cite{Chow:2007} Let $ (\mathcal{M},\mathbf{g}(t)) $ be a solution of RF that exists up to a maximal time $ T \leq \infty $. We say a~solution to the RF develops a:
\begin{enumerate}
\item {\bf Type-I singularity} at a maximal time $ T < \infty $ if
\begin{equation}
\sup_{t \in [0,T)} (T - t)\max\{|\mathbf{Rm}(x,t)| : \: x \in \mathcal{M}\} < \infty; \nonumber
\end{equation}

\item {\bf Type-IIa singularity} at a maximal time $ T < \infty $ if:
\begin{equation}
\sup_{t \in [0,T)} (T - t)\max\{|\mathbf{Rm}(x,t)| : \: x \in \mathcal{M}\} = \infty; \nonumber
\end{equation}

\item {\bf Type-IIb singularity} at a maximal time $ T < \infty $ if:
\begin{equation}
\sup_{t \in [0,\infty)} \max\{|\mathbf{Rm}(x,t)| : \: x \in \mathcal{M}\} = \infty; \nonumber
\end{equation}

\item {\bf Type-III singularity} at a maximal time $ T < \infty $ if:
\begin{equation}
\sup_{t \in [0,\infty)} \max\{|\mathbf{Rm}(x,t)| : \: x \in \mathcal{M}\} < \infty. \nonumber
\end{equation}
\end{enumerate}
\end{defn}

The dumbbells investigated in this work will form nondegenerate ``neckpinch'' singularities. An~example of a nondegenerate neckpinch was established in \cite{AK}. We review the notation and some necessary definitions \cite{AIK:2011,AIK:2015}.

\begin{defn}
\label{def:blowup}
A sequence $ \{(x_{j},t_{j})\}_{j=0}^{\infty} $ of points and times in an RF solution is called a {{blow-up sequence at time}} $ T $ if $ t_{j} \rightarrow T $ and if $ |\mathbf{Rm}(x_{j},t_{j})| \rightarrow \infty $ as $ j \rightarrow \infty $.
\end{defn}

Such a sequence has a corresponding pointed singularity model if the sequence of parabolic dilation metrics $ \mathbf{g}_{j}(x,t) := |\mathbf{Rm}(x_{j},t_{j})|\mathbf{g}(x,t_{j} + |\mathbf{Rm}(x_{j},t_{j})|^{-1}t) $ has a complete smooth limit.

\begin{defn}
\label{def:neckpinch}
Given an RF solution $ (\mathcal{M}, \mathbf{g}(t)) $, a {{neckpinch singularity}} develops at a time $ T $ if there is some blow-up sequence at $ T $ whose corresponding pointed singularity model exists and is given by the self-similarly shrinking Ricci soliton on the cylinder $ \mathbb{R} \times S^{n} $.
\end{defn}

\section{Methodology}
\label{sec:Methodology}
\vspace{-6pt}

\subsection{Models}

We investigate topological signals of singularity formation for a dimpled sphere on $ S^{2} $ and dumbbells on $ S^{3} $ in this work. The types of dumbbells include a symmetric dumbbell and a~dimpled dumbbell (``genetically-modified peanut'' \cite{AK}). 
Singularity formation in such examples is well-understood \cite{AIK:2011,AIK:2015,AK,CK,GZ}.

The ``dimpled'' sphere is modeled by a metric whose radial function depends on the polar angle $ \theta $ and time and is given by:
\begin{equation}
\mathbf{g}(t) = r^{2}(\theta,t)\mathbf{g}_{can}
\end{equation}

\noindent where $ \mathbf{g}_{can} = d\theta^{2} + \sin^{2}\theta d\phi^{2} $ is the canonical metric on a 2-sphere. The angular measures are usually $ -\pi/2 < \theta < \pi/2 $ (polar) and $ 0 \leq \phi < 2\pi $ (azimuthal). For our algorithm, it is convenient to modify the polar angle so that $ 0 < \theta < \pi $; this is the angle over which we construct the initial radial profile~$ r(\theta,0) $.

By the uniformization theorem \cite{Hamilton:1982}, any 2-geometry will evolve under RF to a constant curvature sphere, plane or hyperboloid. Regarding singularity formation, this will lead to the collapse of a~sphere to a round point, and the dimpled sphere will Ricci flow to a sphere of constant curvature before collapsing to a round point. The shrinking round sphere is a Type-I singularity.

The dumbbells are constructed by puncturing $ S^{3} $ at the poles $ \{N,S\} $ and identifying $ S^{3} \setminus \{N,S\} $ with $ (-c,c) \times S^{2} $ where $ c $ is a constant usually taken as 1. This identification is for convenience to avoid working in multiple patches \cite{AIK:2011,AIK:2015}. Letting $ x $ denote the coordinate on $ (-c,c) $ and $ \mathbf{g}_{can} $ denote the canonical unit sphere metric on $ S^{2} $, an arbitrary family $ \mathbf{g}(t) $ of smooth $ SO(3) $-invariant metrics on $ S^{3} $ may be written in geodesic polar coordinates as:
\begin{equation}
\label{eq:geodesicpolarmetric}
\mathbf{g}(t) = \varphi^{2}(x,t)dx^{2} + \psi^{2}(x,t)\mathbf{g}_{can}.
\end{equation}

The function $ \psi^{2}(x,t) $ is a ``radial'' function for the profile of the dumbbell with one dimension suppressed. This manifests as corseting of the dumbbell.

An alternative representation of this metric, the warped product metric, is to introduce a~geometric coordinate which normalizes the metric in the variable associated with the interval and leaves only a radial function present \cite{AK}. The distance from the ``equator'' $ \{0\} \times S^{n} $ is given by:
\begin{equation}
s(x,t) = \int_{0}^{x}\varphi(x',t)dx'.
\end{equation}

It follows that $ \partial s/\partial x = \varphi(x,t) $, and this changes the metric to:
\begin{equation}
\label{eq:normalizedmetric}
\mathbf{g}(t) = ds^{2} + \psi^{2}(s,t)\mathbf{g}_{can}.
\end{equation}

We formulate an initial-value and boundary-value problem for each of the dumbbells; 
this~is achieved conveniently by considering both $ \varphi(x,t) $ and $ \psi(x,t) $, so we use (\ref{eq:geodesicpolarmetric}). 
The initial and boundary conditions for the radial function $ \psi(x,t) $ of the dumbbells depend on the model. We~use Neumann conditions on the boundaries to satisfy the requirements for the values of the derivatives on approaching the poles \cite{AK},
\begin{equation}
\label{eq:derivreq}
\lim_{x\rightarrow\pm c} \frac{\partial\psi}{\partial s} = \mp 1,
\end{equation}

\noindent which can be expressed in terms of $ x $ since $ (\partial x/\partial s)|_{t=0} = 1/\varphi(x,0) $. The boundary conditions at the poles are necessary such that the metric on the dumbbell extends smoothly to a metric on $ S^{3} $.

The dumbbells, under RF, will manifest neckpinch singularities. These come in two types: {nondegenerate} and {degenerate}.

\begin{defn}
\label{def:nondeg}
A neckpinch singularity is {{nondegenerate}} if every pointed singularity model of any blow-up sequence corresponding to $ T $ is a cylindrical solution. The following are basic assumptions in \cite{AIK:2011,AIK:2015} for such singularity formation in $ SO(n+1) $-invariant solutions of RF with an initial set of data of the form (\ref{eq:normalizedmetric}): 

\begin{enumerate}
\item The sectional curvature $ L $ of planes tangent to each sphere $ \{s\} \times S^{n} $ is positive.
\item The Ricci curvature $ \mathbf{Rc} = nK ds^{2} + [K + (n-1)L]\psi^{2}\mathbf{g}_{can} $ (where $ K $ is the sectional curvature of a plane orthogonal to $ \{s\} \times S^{n} $) is positive on each polar cap.
\item The scalar curvature $ R = 2nK + n(n-1)L $ is positive everywhere.
\item The metric has at least one neck and is ``sufficiently pinched'' in the sense that the value of the radial function $ \psi $ at the smallest neck is sufficiently small relative to its value at either adjacent bump. 
\item The metric is reflection symmetric, and the smallest neck is at $ x = 0 $. 
\end{enumerate}
\end{defn}

Here,
\begin{eqnarray}
L &=& \frac{1 - \psi_{s}^{2}}{\psi^{2}} = \frac{1 - (\psi_{x}/\varphi)^{2}}{\psi^{2}}, \:\:\: \textrm{and} \nonumber\\
K &=& -\frac{\psi_{ss}}{\psi} = -\frac{(-\varphi_{x}\psi_{x}/\varphi^{2} + \psi_{xx}/\varphi)}{\varphi\psi}.
\end{eqnarray}

Nondegenerate neckpinches are Type-I singularities.

An obvious consequence of the positivity of the tangential sectional curvature $ L $ is that \mbox{$ |\psi_{s}| \leq 1$}. This makes construction of suitable models more challenging and places constraints on any randomized~construction.

\subsection{Persistence Computations 
}
\label{perscomp}

Consider a triangulated manifold $\mathcal{K}$ of dimension $n$ such that each $d$-dimensional simplex $\sigma$ has been assigned a real number $\iota(\sigma)$. From this data, one can construct a filtration $\{\mathcal{K}_p \mid p \geq 0\}$ {compatible} with $\iota$ such that for each $\sigma$ in $\mathcal{K}$ of dimension $d$, 
\[
\inf ~ \{p \in \mathbb{R} \mid \sigma \in \mathcal{K}_p\} = \iota(\sigma).
\]

That is, the smallest value of $p$ for which a given $d$-simplex $\sigma$ lies in $\mathcal{K}_p$ is required to equal $\iota(\sigma)$. Various ad hoc ways exist to extend $\iota$-values to simplices of dimensions other than $d$ so that a~compatible filtration may be imposed on $\mathcal{K}$. A more principled route is via the {star filtration} which has been used, for instance, in the persistent homological analysis of image data \cite{Robins}.

\begin{defn}
\label{def:starfiltration}
Given a triangulated $n$-manifold $\mathcal{K}$ and a real-valued function $\iota$ defined on the $d$-simplices for some $0 \leq d \leq n$, the {{star filtration}} along $\iota$ is defined via the following containment relation for each $p \in \mathbb{R}$ and each simplex $\sigma \in \mathcal{K}$: we have $\sigma \in \mathcal{K}_p$ if and only if one of the following conditions holds:
\begin{itemize}
\item $\dim \sigma < d$ and there is some $d$-dimensional co-face $\tau \succ \sigma$ with $\iota(\tau) \geq p$, or
\item $\dim \sigma = d$ and $\iota(\sigma) \leq p$, or
\item $\dim \sigma > d$ and there is some $d$-dimensional face $\tau \prec \sigma$ with $\iota(\tau) \leq p$.
\end{itemize}
\end{defn}
The extreme cases $d = 0$ and $d=n$ are called upper and lower star filtrations. To construct a star filtration, consider a simplex in $\mathcal{K}$. If its dimension equals $d$, then its filtration index equals its $\iota$-value. If its dimension exceeds $d$, then it inherits the largest $\iota$-value encountered among its $d$-dimensional faces. Finally, if its dimension is smaller than $d$, then it inherits the smallest $\iota$ value encountered among all $d$-dimensional simplices in $\mathcal{K}$ which contain this simplex as a face. Note that in order for such a construction to be well-defined, each simplex of dimension smaller than $d$ must have at least one $d$-dimensional co-face. However, since we work entirely with triangulated manifolds, this condition is automatically satisfied.

The following 
describes our methodology for analyzing RF with PH. We begin with a triangulated surface $\mathcal{K}^0$ evolving via RF with the resulting complexes being labeled $\mathcal{K}^t$ for some discrete indices~$t \geq~0$:
\begin{enumerate}
\item for each $\mathcal{K}^t$, compute curvature values assigned to all vertices,
\item construct the upper star filtration along these values,
\item produce the corresponding PDs in dimensions $0$ and $1$,
\item instead of birth-death pairs $(b,d)$, restrict attention to the
differences $d-b$, called {persistence intervals} or {lifespans}.
\end{enumerate}

The superscript, associated to discrete time, distinguishes from the subscript, associated with the filtration. {Perseus} \cite{Pers} generates a reduced complex $ \mathcal{K'}^t $ out of the input simplicial complex $\mathcal{K}^t$. The output files depend on the dimension of $ \mathcal{K'}^t $, which is less than or equal to the dimension of $\mathcal{K}^t$. The output files come in two forms at each discrete time step $t \geq 0$. One form consists of tables containing persistence intervals for each dimension of $ \mathcal{K'}^t $. The other form is a table of the Betti numbers $ \beta_{i} $, $ i = 0,...,n = \dim(\mathcal{K'}^t) $, 
at each value of the filtration parameter (here, scalar curvature). 
These numbers provide a count of the number of connected components ($ n = 0 $) or holes of increasing dimensions ($ n > 0 $). Then, we plot the lifespans against 
discretized values of 
$t \geq 0$. 

In tandem with the methodology above, we consider PDs--multisets of points each identified with the birth and death value of a lifespan, $ (b,d) $, in the extended Euclidean plane $ \widetilde{\mathbb{R}}^{2} $ \cite{EH:2008}---associated with these tables and compute several distance measures between them: bottleneck ($ d_{B} $), Wasserstein-1 ($ d_{W^{1}} $) and Wasserstein-2 ($ d_{W^{2}} $). The {bottleneck distance} is defined as:
\newpage
\vspace{12pt}

\begin{equation}
d_{B} = \inf_{\gamma}\sup_{p\in PD_{i}}||p - \gamma(p)||_{\infty}
\end{equation}

\noindent where $ \gamma: PD_{i} \rightarrow PD_{j} $ is a distance-minimizing bijection between PDs, here indexed finitely, and~\mbox{$ \displaystyle ||p - \gamma(p)||_{\infty} = \max\{|p_{x} - \gamma(p)_{x}|,|p_{y} - \gamma(p)_{y}|\} $} is the $ L_{\infty} ${-norm} on $ \widetilde{\mathbb{R}}^{2} $. For the $ L_{\infty} $-norm, distances between infinite coordinates are defined to be zero. The {Wasserstein-q} distance is: 
\begin{equation}
d_{W^{q}} = \left(\sum_{p\in PD_{i}}||p - \gamma(p)||_{\infty}^{q}\right)^{1/q}.
\end{equation}

Note that $ \displaystyle \lim_{q \rightarrow \infty} d_{W^{q}} = d_{B} $. For our system of PDs of finite cardinality (compared with PDs of continuous functions), in the bottleneck distance, the infimum is a minimum, and the supremum is a~maximum.

The bottleneck distance $ d_{B} $~measures the largest difference between PDs, while the Wasserstein -1 and -2 distances $ d_{W^{1}} $ and $ d_{W^{2}} $ measure all of the differences (with $ d_{W^{2}} $ less sensitive to small changes than $ d_{W^{1}} $) \cite{KGKM:2014}.  If $ d_{W^{1}} $ is significantly larger than $ d_{B} $ and $ d_{W^{2}} $, then many small features are responsible for changes in the geometry between PDs. If the distance measures are close, then this indicates a single dominant feature is responsible for changes in the geometry between PDs. The Betti number spectra, tables measuring lengths of persistence intervals, and distances between persistence diagrams together provide a topological signal associated to a collection of data.

For the models, we investigate how the distance measures change over the time evolution of the system. This approach allows us to acquire persistent homological characterizations of the spaces---hence topological signals of the formation and types of singularities---as they evolve via RF.

\subsection{Data Generation and Preparation Algorithm}
\label{data}

By numerically solving the Ricci flow equations, we compute the scalar curvature. The choices of initial and boundary conditions depend on the particular model.

The initial radial profile for the dimpled sphere with topology $ S^{2} $ is obtained by constructing a randomized table of values for a collection of angles $ \theta $. To solve RF, we evolve over the angle $ \theta $ and time. The interval for $ \theta $ is $ [\epsilon, \pi - \epsilon] $ with $ \epsilon = 10^{-3} $. This choice of interval allows for more stable evolution by reducing potential complications associated with boundary effects.

To obtain a numerical solution in the evolution of the dumbbells, all of which have topology $ (-c,c) \times S^{2} $, we approximate the topology of the dumbbells as $ [-c_{Boundary},c_{Boundary}] \times S^{2} $. The~endpoints of the closed interval are $ c_{Boundary} = c - \varepsilon $ where $ \varepsilon = \mu\epsilon $, and $ \mu $ is a scale factor. The~endpoints $ \pm c_{Boundary} $ are, variously, $ \pm \mu(\pi/2 - \epsilon) $ (symmetric dumbbell) and $ \pm \mu(\pi - \epsilon) $ (dimpled~dumbbell).

For all dumbbells, the function $ \varphi(x,t) $ is given initially by $ \varphi(x,0) = 1 $ and at the boundaries by $ \varphi(c_{Boundary},t) = 1 + t $. This choice initially satisfies the boundary conditions (\ref{eq:derivreq}) for $ \psi(x,t) $ on the poles since $ (\partial x/\partial s)|_{t=0} = 1/\varphi(x,0) = 1 $.

To address properly the work of \cite{GZ} (used in Section \ref{sec:Results}), we need for the initial profile a function, depending on a parameter $ \alpha \in [0,1] $, that satisfies the derivative requirements on the poles for the smooth extension of the metric to $ S^{3} $; we impose this on $ g_{\alpha}(0) $. Such a function for the initial profile is given by:
\begin{equation}
\label{eq:alphainterpolation}
\psi(x,0) = \frac{-\mu a\left(\frac{x}{\mu} + \frac{\pi}{2}\right)\left(\left(\frac{x}{\mu} - (1 - \alpha)L\right)^{2} + h(\alpha,k)\right)\left(\frac{x}{\mu} - \left(\left(\frac{\pi}{2} - L\right)\alpha + L\right)\right)}{(b - \frac{x}{\mu})^{(1 - \alpha)}}.
\end{equation}

We have a freely-specifiable set of parameters $ \{\alpha, L, k, \mu\} $. The factor $ k $ appears in the definition of $ h(\alpha,k) = -k\alpha + 2k $; this function controls the neck thickness. The factor $ (x/\mu - (1 - \alpha)L)^{2} + h(\alpha,k) $ controls the position of the neck. We select $ \mu = 100 $, and $ \epsilon = 10^{-3} $. This choice of $ \mu $ allows for longer time evolution and delay of stiffness due to incrementation limitations. 
The constants $ a $ and $ b $ depend on the choice of these parameters and are determined by imposing the condition of the derivatives on the poles. 
The left pole is fixed; 
the right endpoint is a function of $ \alpha $ and $ L $, $ \mu\left(\left(\frac{\pi}{2} - L\right)\alpha + L\right) $.

To consider nondegenerate neckpinches satisfying the conditions of \cite{AIK:2011,AIK:2015}, we consider a rotationally and reflection symmetric (about the ``equator'') dumbbell. This is obtained from (\ref{eq:alphainterpolation}) with $ \alpha = 1 $ and parameters $ L = 1 $ and $ k = 1/25 $ providing for sufficient pinching. Furthermore, this model satisfies the conditions of Definition \ref{def:nondeg}. The initial radial profile is given by:
\begin{equation}
\psi(x,0) = -12.6948\left(\frac{x}{100} + \frac{\pi}{2}\right)\left(\left(\frac{x}{100}\right)^{2} + (1/25)\right)\left(\frac{x}{100} - \frac{\pi}{2}\right).
\end{equation}

Interesting topological features occur when the filtration parameter or function has many critical points (e.g., the height function \cite{EH:2008}). A way to achieve such a function is to consider the dimpled dumbbell or ``genetically modified peanut'' \cite{AK}. This is an object with multiple necks and bumps. For~the dimpled dumbbell, we construct an interpolation function over the interval $ [-100\pi,100\pi] $. This interval is larger than for the other dumbbells and is chosen for more stable evolution (though we still come in by 0.1 unit on each pole for the RF computation) and to ease placement of the necks while respecting the derivative conditions for dumbbells.

\begin{figure}[H]
\centering
\includegraphics[height=3.75in,width=5.75in]{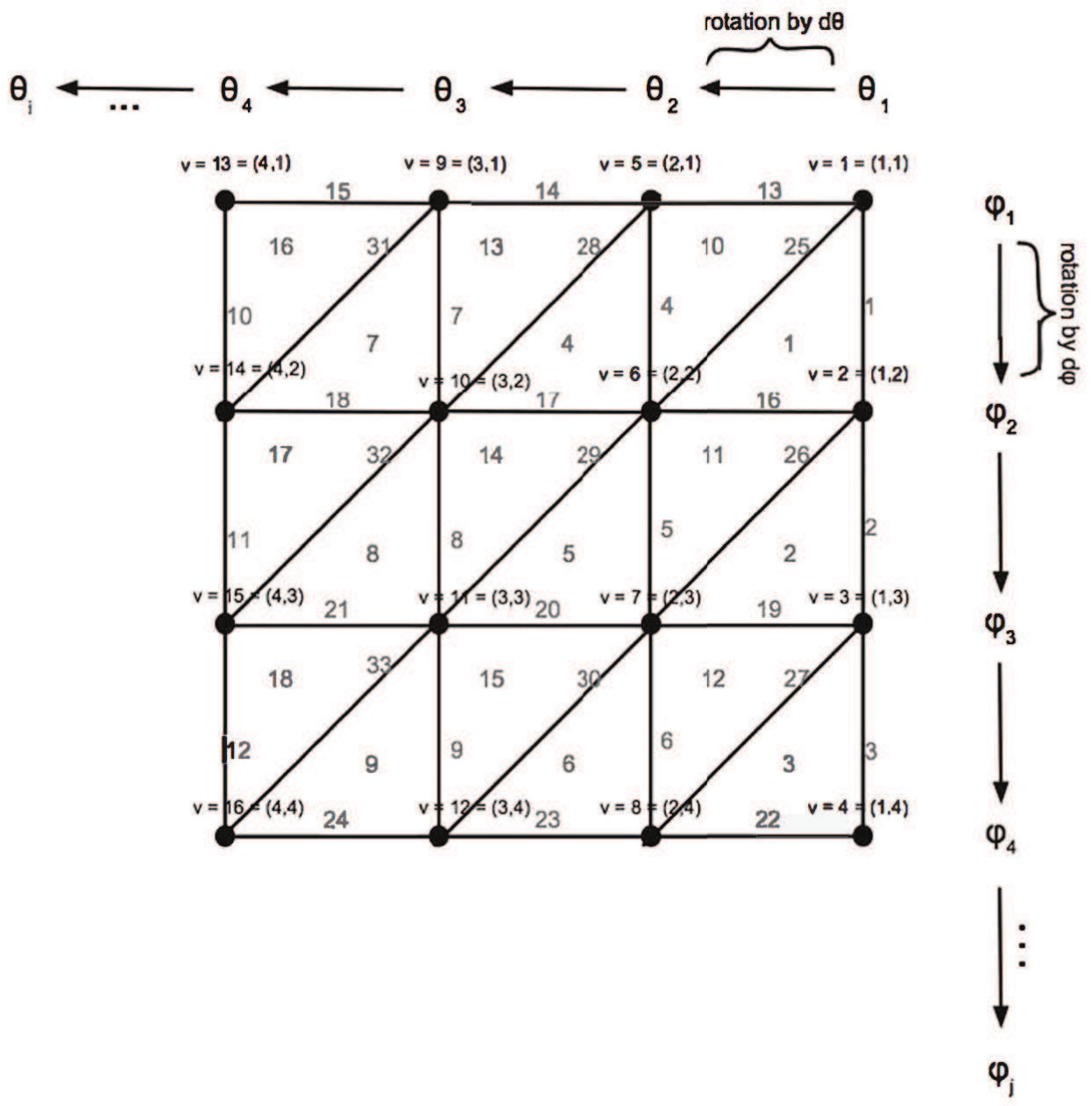}
\caption{Example triangulation where the index $ i $ corresponds with the polar angle $ \theta $, and the index $ j $ corresponds with the azimuthal angle $ \phi $. Vertices are indexed as $ (i,j) $ or as the single digit $ v $, given by $ v = j + (i-1)N\phi $. 
The numbers of points sampled along the directions are given by the notation $ N\theta $ or $ N\phi $. The increments of angles are computed by $ \triangle\theta = \frac{\pi - 2\varepsilon}{N\theta - 1} $ and $ \triangle\phi = \frac{2\pi}{N\phi - 1} $.}
\label{fig:TriDiagram}
\end{figure}

To prepare the data for the PH analysis, we triangulate the manifold in question (Figure~\ref{fig:TriDiagram}). Other~suitable triangulations of a manifold for numerical evolution of Ricci flow exist (cf. \cite{Miller:2014crf}); we~have chosen this one for convenience. For each vertex chosen, we obtain the scalar curvature. To~compute the scalar curvature along a given edge of the triangulation, we average the scalar curvatures of the endpoints of the edge. This information is then used to determine scalar curvature values of all edges of the triangulation for each time index.

{Perseus} implements this data to compute persistence intervals. To do so, the vertices, edges, and~faces must be indexed such that boundaries of the edges (vertices) correspond with the appropriate edges, and boundaries of faces (edges) correspond with the appropriate faces. \mbox{Then, this information} is input into matrices that {Perseus} can efficiently reduce.

For all of our models, we are able to implement the same triangulation. This follows from the symmetry of the dumbbell: in computing the components of the Ricci tensor, the $ \theta\theta $ and $ \phi\phi $ components are redundant so that, out of three, only the $ xx $ and $ \theta\theta $ components need be computed to determine the scalar curvature. The suppression of one of the dimensions of $ S^{2} $ thus allows the use of a similar triangulation to the dimpled sphere. In the triangulation, the modifications from the dimpled sphere to the dumbbells are $ \theta \rightarrow x $ and:
\begin{eqnarray}
\triangle x &=& 100\frac{\pi - 2\varepsilon}{Nx - 1}.
\end{eqnarray}


Furthermore, because the radial function for the dimpled sphere possesses no dependence on the angle $ \phi $, one can consider an arbitrarily coarse or fine incrementation for $ \phi $; the results from {Perseus} will remain the same. The same is true regarding the dumbbells whose radial functions have no dependence on the angle $ \theta $. Hence, in our models, we will identify them according to the number of points used to generate the partitions (e.g., $ N\theta $ and $ Nx $ for the sphere and dumbbells, respectively).

We use a common scheme for generating a discrete sample of time indices for each of the models. This involves counting the first ten time indices in units of one, the next ten in units of ten, then units of a hundred, then units of a thousand, and then units of ten-thousand. After this, we shorten the increments variously, depending on the model, so that we end by looking at changes again in units of~one.

Following this, we compute the pairwise distances between adjacent PDs using the $ L_{\infty} $-norm. The approach taken here is to take a point within, for example, $ PD_{i} $ and to compute its distance (in~the~$ L_{\infty} $-norm) from each of the points in $ PD_{i+1} $. We repeat this for all of the points in $ PD_{i} $. If there is a difference in cardinalities of the PDs, we can take points from the diagonal of the set of lesser cardinality to make these equal. Computationally, this produces an $ n \times n $ square matrix where $ n $ is the cardinality of the (larger of the two) PDs. When we encounter such a situation, we estimate the point that yields a minimizing bijection for the distance measures.

In the tables for the connected components, one lifespan will have an infinite death value corresponding with the first feature born and last feature to die. In a PD, each point is identified with the birth and death value of a lifespan, $ (b,d) $. As stated above, using the $ L_{\infty} $-norm, the distance between infinite coordinates is defined to be zero. The distance from a point with two finite coordinates and a point with a finite coordinate and an infinite coordinate, in the $ L_{\infty} $-norm, will be infinite (since~the~difference between infinity and a finite number is again infinity). Thus, the matrix of $ L_{\infty} $ distances of adjacent PDs for the connected components will have a column and row of entries such that all but one will be infinity. This is not a consideration for the adjacent PDs for holes since no lifespan is infinite. {Perseus} generates tables for the holes which have multiple copies of the same lifespans; we delete these duplicates to construct a PD at each time index where each distinct lifespan corresponds with a single point in the PD.

Then, we compute the bottleneck, Wasserstein-1, and Wasserstein-2 distances. We collect together the bijections of the PDs by taking permutations of the square matrix of $ L_{\infty} $ distances. This allows us to consider the bijections as diagonals of permuted matrices, enabling convenient computation of searches of maximum distances (along the diagonal in the case of the bottleneck distance) and in computing sums (traces of a matrix for the Wasserstein distances). For the connected components, that~the permutations can be restricted to the permutations of the $ (n-1) \times (n-1) $ matrix corresponding to the $ L_{\infty} $ distances between points with only finite coordinates. This is because the distance between points with only finite coordinates and a point with an infinite coordinate will be infinity, which would inherently be the maximum distance. However, the bottleneck distance looks for the minimum across all of the maximal distances, so none of these infinities will be a minimum. The Wasserstein distances follow~similarly.

\section{Results}
\label{sec:Results}

First, we plot the initial and final radial profiles for the models. For the models, we consider two triangulations: one coarser and one finer. Then, we plot the associated lifespans and investigate the outputs from the different data tables generated by {Perseus}. We refer to these plots and tables in terms of Betti number so that we have plots and tables associated to $ \beta_{0} $ (number of connected components) (we omit similar plots for $ \beta_{1} $ (the number of tunnels) since they do not provide more or distinct information). The symmetry of our models means that we have only dimensions zero and one.

Furthermore, we compute the pairwise bottleneck ($ d_{B} $), Wasserstein-1 ($ d_{W^{1}} $) and Wasserstein-2 ($ d_{W^{2}} $) distances between adjacent PDs. We plot these and the average scalar curvature at each time index. Then, we compute the ratios $ d_{B}/d_{W^{1}} $ and $ d_{B}/d_{W^{2}} $, both always less than or equal to one due to $ d_{B} $ being the smallest, to obtain a quantification and characterization of topological signal as we move across time indices.

\subsection{Dimpled Sphere}
\label{ds}

\begin{figure}[H]
\centering
\begin{tabular}{cc}
 \includegraphics[width=2.5in]{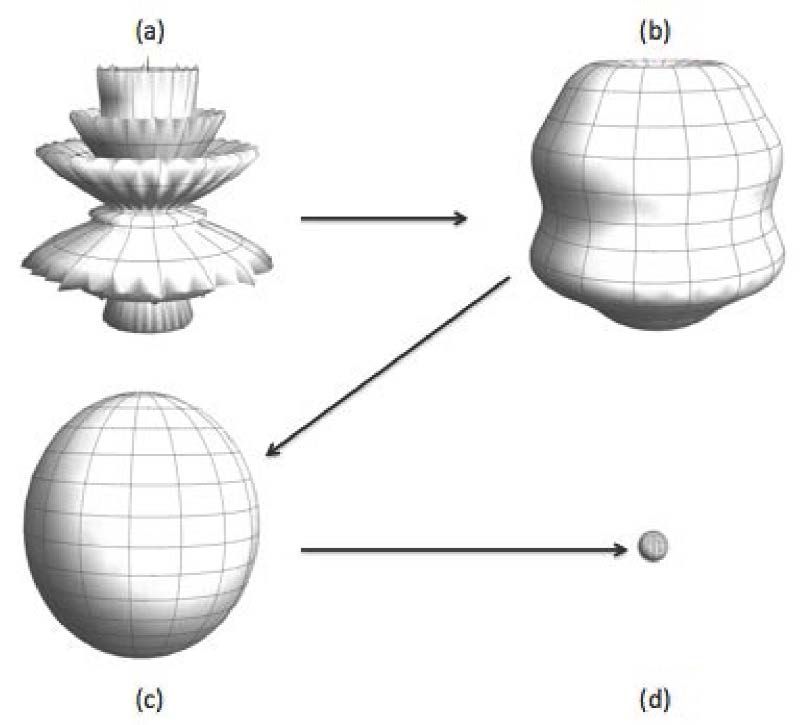} &
 \includegraphics[width=2.5in]{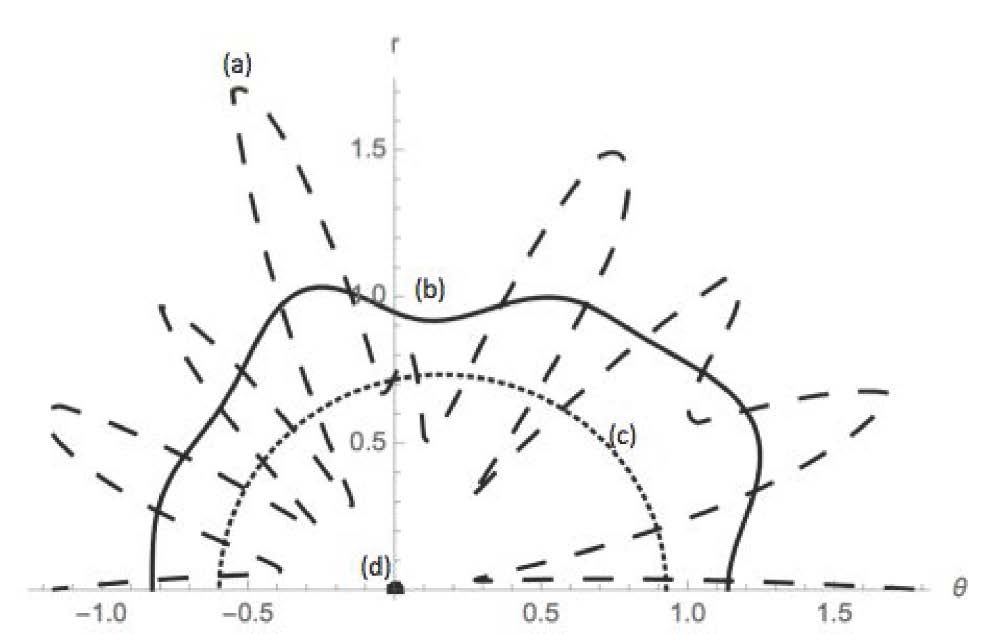}
 \end{tabular}
\caption{This illustration is of the dimpled sphere and the corresponding radial profiles
 $ r(\theta,t) $ at time ({\textbf{a}}) $ t = 0 $; ({\textbf{b}}) $ t = 0.03 $; ({\textbf{c}}) $ t = 0.3 $; and ({\textbf{d}}) $ t = 0.57105 $.}
\label{fig:2DRotSolid}
\end{figure}

The rotational solid in Figure~\ref{fig:2DRotSolid} has the topology of $ S^{2} $ and the appearance of a ``dimpled'' sphere. The dimpling is obtained by using a random number generator to modulate the initial radius of the solid so that $ r(\theta,0) $ is an interpolation function built over a table of angles and radii. Ricci flow uniformizes the curvature, transitioning the dimpled sphere to a more conventional-looking sphere, which then shrinks to a round point. The triangulations are taken over $ N\theta = 15 $ and $ N\theta = 25 $ grids for 78 time indices. While finer triangulations can be considered, they result in significantly longer computational times to obtain the various distance measures.

\begin{figure}[H]
\begin{center}
 \begin{tabular}{cc}
 \includegraphics[width=2.5in]{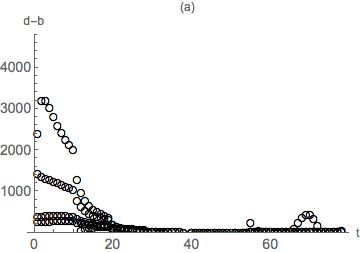} &
 \includegraphics[width=2.5in]{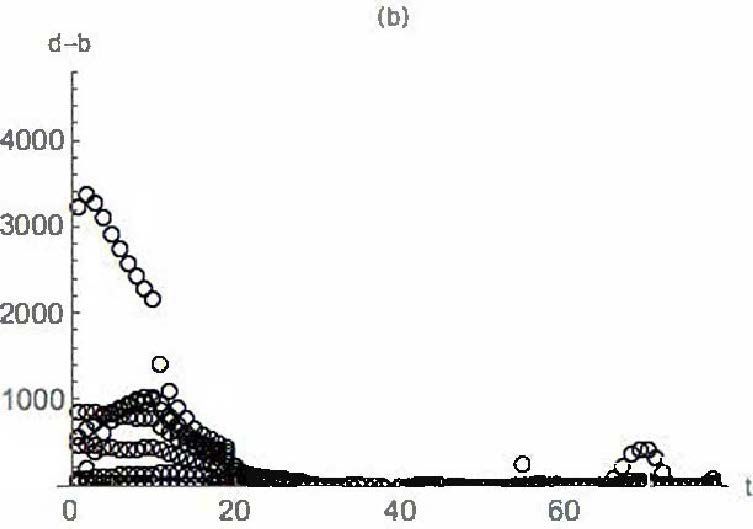}
 \end{tabular}
\caption{Lifespans ({d}--{b}) for $ \beta_{0} $ computed for ({\textbf{a}}) $ N\theta = 15 $
 and ({\textbf{b}}) $ N\theta = 25 $ triangulations.}
\label{fig:2DLifespans}
\end{center}
\end{figure}

\begin{figure}[H]
\begin{center}
 \begin{tabular}{cc}
 \includegraphics[width=2.5in]{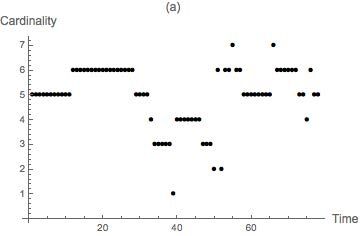} &
 \includegraphics[width=2.5in]{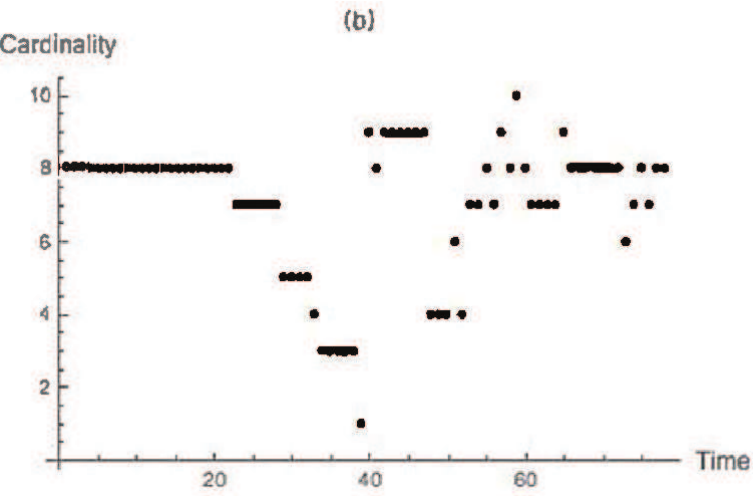}
 \end{tabular}
\caption{Cardinalities (number of points in PDs) at each time index for $ \beta_{0} $ computed for ({\textbf{a}}) $ N\theta = 15 $ and ({\textbf{b}}) $ N\theta = 25 $.
}
\label{fig:CardinalityDS}
\end{center}
\end{figure}

Figure \ref{fig:2DLifespans} features the lifespans $ d - b $ (death-birth) associated to the Betti number $ \beta_{0} $ versus time index $ t $ for the triangulations. For all time indices, the tables for connected components have a single infinite lifespan. Initially, the tables for connected components and holes have multiple finite lifespans. We plot these in Figure \ref{fig:CardinalityDS}. 

The numbers of lifespans decrease until, at time $ t=0.3 $ (time index 39 for the connected components), the $ N\theta = 15 $ and $ N\theta = 25 $ triangulations have no finite lifespans and only an infinite lifespan in their table for connected components. Therefore, they do not have any lifespans for their tables of holes. For a $ N\theta = 50 $ triangulation (not plotted), the table of connected components has two~finite lifespans, of length essentially zero. The table of holes at this time index also has two lifespans: one of length essentially zero and another of negligible length. In comparing this with the lifespans measured at the first time index, we find that the lifespans have decreased by three orders of magnitude for $ \beta_{0} $ and two orders of magnitude for $ \beta_{1} $. These findings suggest that at the scale of resolution for the coarser triangulation, the reduced complex does not generate any noise in the process of increasing the filtration parameter. Geometrically, the interpretation is that the scalar curvature throughout the solid has become uniform (numerically) and begun its collapse.

\unskip

As the solid collapses to a round point, the number of lifespans increases until, at the end of the evolution (when our RF algorithm stiffens), we have four finite lifespans for the connected components and holes for $ N\theta = 15 $ and seven finite lifespans for the connected components, but six for the holes for $ N\theta = 25 $. From experiments (not shown) with a typical sphere whose radius is the same for all values of $ \theta $, we encounter a similar situation initially where the number of finite lifespans is zero for the connected components. As with the dimpled sphere, this is consistent with the sphere having uniform curvature. Also as with the dimpled sphere, the number of lifespans grows as collapse to a round point occurs. The interpretation is that the numerical code evolves the data points with slight nonuniformity (with respect to scalar curvature), which is in some sense noise.

In the tables at each time index, we have one connected component and zero holes for the lowest and highest maximum birth scalar curvature with varying numbers of connected components and holes for the scalar curvatures between these values. These varying numbers indicate noise related to the process of filtration. Furthermore, at each time index approaching 30,000, the difference between the lowest and highest maximum birth scalar curvatures becomes smaller; this is an indication that the scalar curvature of the object has become uniform. Under RF, it is expected that the dimpling smooths so that the object attains a uniform scalar curvature, then collapses.

\begin{figure}[H]
\begin{center}
 \begin{tabular}{cc}
 \includegraphics[width=2.8in]{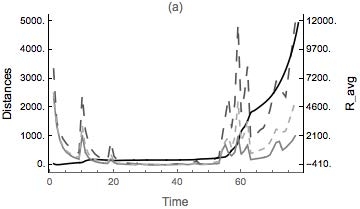} &
 \includegraphics[width=2.8in]{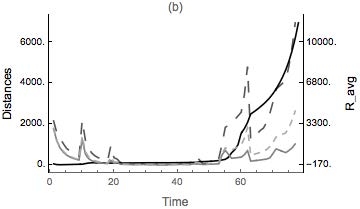}
 \end{tabular}
\caption{$ \beta_{0} $: interpolation functions of average scalar curvature (denoted by R\_avg, solid black line), bottleneck (gray), Wasserstein-1 (dark gray, dashed), and Wasserstein-2 (light gray, dashed) distances for ({\textbf{a}}) $ N\theta = 15 $ and ({\textbf{b}}) $ N\theta = 25 $.}
\label{DSDists}
\end{center}
\end{figure}

We examine the three distance measures in Figure \ref{DSDists} ($ N\theta = 15 $ and $ N\theta = 25 $). These are plotted along with the average scalar curvature as measured at each sampled time. The average scalar curvature is computed with respect to the triangulations, sampling the points accordingly. Approaching the time index where the scalar curvature becomes uniform across the object, the~distances all decrease; after this time index, the distances increase again. The average scalar curvature does not change until after uniformization, where it increases significantly. Likewise, $ d_{W^{1}} $ becomes particularly large relative to $ d_{W^{2}} $ and $ d_{B} $, an indication that many small changes are occurring with the geometry. 
Given each time index, the adjacent PDs represented by that index have the common feature that $ d_{W^{2}} $ is closer to $ d_{B} $ than it is to $ d_{W^{1}} $. 

\begin{figure}[H]
\begin{center}
 \begin{tabular}{cc}
 \includegraphics[width=2.8in]{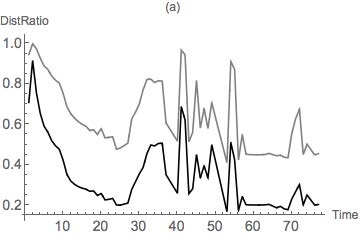} &
 \includegraphics[width=2.8in]{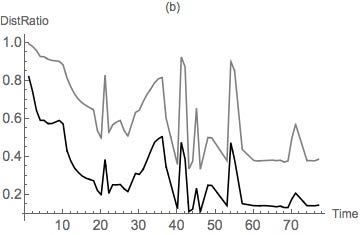}
 \end{tabular}
\caption{$ \beta_{0} $: Interpolation functions of ratios of bottleneck distance to Wasserstein-1 distance (black dots) and Wasserstein-2 distance (gray dots) for ({\textbf{a}}) $ N\theta = 15 $ and ({\textbf{b}}) $ N\theta = 25 $.}
\label{DSDistRat}
\end{center} 
\end{figure}

To quantify more precisely the similarities and differences between the distance measures, and to compensate for scaling, we compute the ratios of the distance measures $ d_{B}/d_{W^{1}} $ and $ d_{B}/d_{W^{2}} $ for the connected components for the adjacent PDs (Figure \ref{DSDistRat}) for both triangulations (the same calculations and plots for holes reveals no more information, so we do not display them). Initially, we find the ratios in close proximity, especially for the coarser triangulation. This is consistent with the sharp dimpling present in the beginning as such dimpling consists of large changes in the geometry. For both the connected components and the holes, the distance ratios fluctuate but trend downward over the entire evolution. Spikes occur in relation to the changes of time incrementation, which results in larger differences between persistence diagrams. Near the time of the scalar curvature becoming uniform, the distance ratios reach their lowest points up to that time. After this, they remain low, save for a few spikes. This is an indication that many small changes contribute to changes in the geometry of the sphere as it collapses to a round point. Naive intuition would suggest that collapse to a round point is a single large effect (thinking in terms of a continuous object), but as we are looking at a discrete collection of data, and because PH analyzes the data points, this collapse is interpreted in terms of the behavior of the data points. For the discrete data points to collapse to a round point, the points must tend to a common point. Observe that the lifespans, at the end of the recorded evolution, are increasing but still close in their values. For collapse to a round point, the indication from the topological signal is that the collapse is due to many small features in the data set changing at each time step.

 \unskip

\subsection{Nondegenerate Neckpinch}
\label{nn}
\vspace{-6pt}

\subsubsection{Symmetric Dumbbell}

\begin{figure}[H]
\centering
 \includegraphics[width=5in]{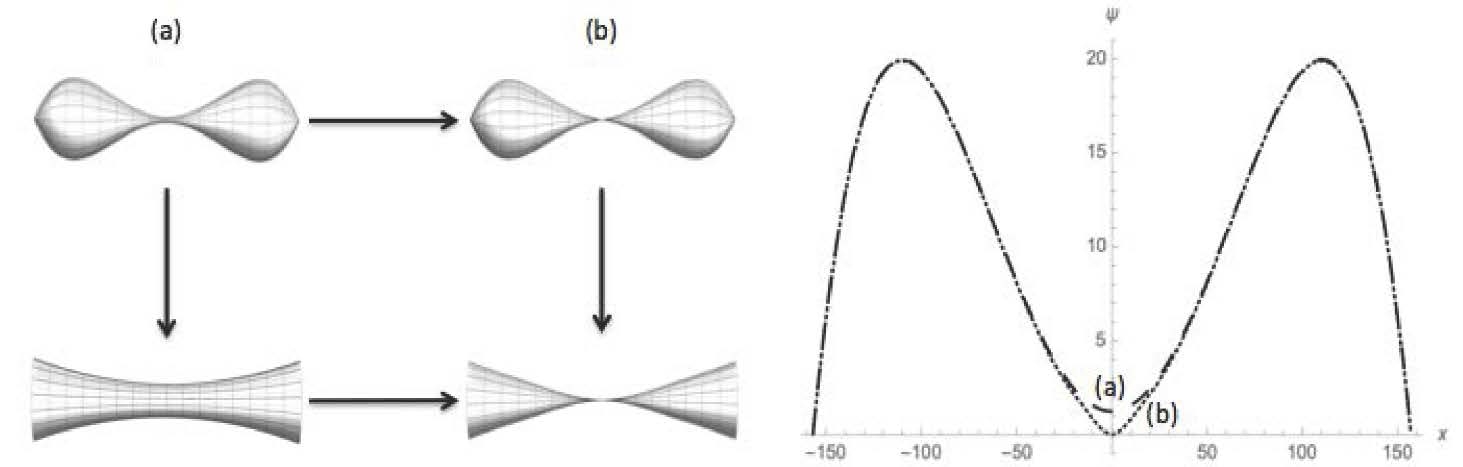}
\caption{This is an illustration of the symmetric dumbbell. The three-dimensional surfaces (above) are the plots of the dumbbell on the interval $ [-50\pi,50\pi] $ and close-ups of the necks at initial time $ t = 0 $ ({\textbf{a}}) and final time $ t = 0.785 $ ({\textbf{b}}). The initial (long dash) and final (short dash) radial profiles are plotted below on the interval $ [-50\pi,50\pi] $.}
\label{fig:SymDumbbell} 
\end{figure}

Figure \ref{fig:SymDumbbell} is a plot of the initial and final surfaces, their corresponding neckpinches, and the initial and final radial profiles. The neck radius (at $ x = 0 $) 
decreases three orders of magnitude over the time evolution. The rest of the plot is essentially static with few changes in the values of $ \psi $. We consider two resolutions for triangulations: $ Nx = 50 $ (coarse-grained) and $ Nx = 100 $ (finer-grained).

For the 
tables of connected components, we have one finite lifespan and one infinite lifespan 
of low scalar curvature for all time indices. 
The presence of only a single finite lifespan indicates that we have only one large feature driving the change in curvature of the dumbbell. The 
tables of holes have multiple lifespans of negligible length. The appearance and disappearance of tunnels indicates that these features are noise related to the filtration process; after we have accounted for all values of curvature, these extraneous features vanish. Over the course of the evolution, for $ Nx = 50 $, the~lifespans increase by two orders of magnitude while for $ Nx = 100 $, they increase by three orders of~magnitude.

\begin{figure}[H]
\centering
 \begin{tabular}{cc}
 \includegraphics[width=2.5in]{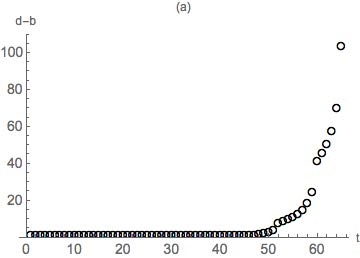} &
 \includegraphics[width=2.5in]{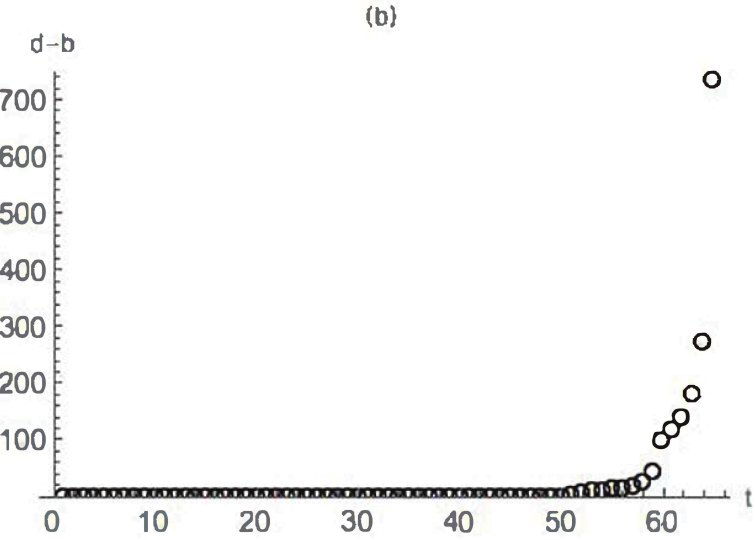}
 \end{tabular}
\caption{Lifespans $ d - b $ for $ \beta_{0} $ computed for ({\textbf{a}}) $ Nx = 50 $ and ({\textbf{b}}) $ Nx = 100 $ triangulations.}
\label{fig:DumbbellPersistence}
\end{figure}

Due to the points sampled in the triangulation, the coarser $ Nx = 50 $ triangulation does not capture the high-curvature value for the pinching which occurs at and near the origin. However, the~finer $ Nx = 100 $ triangulation contains points with higher curvature values. Despite these differences, the plots of the lifespans (Figure \ref{fig:DumbbellPersistence}) for coarser and finer triangulations, respectively, appear almost self-similar.

\begin{figure}[H]
\begin{center}
 \begin{tabular}{cc}
 \includegraphics[width=2.5in]{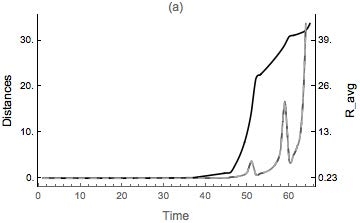} &
 \includegraphics[width=2.5in]{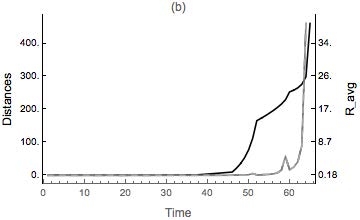}
 \end{tabular}
\caption{$ \beta_{0} $: Interpolation functions of average scalar curvature (denoted by R\_avg, solid black line), bottleneck (gray), Wasserstein-1 (dark gray, dashed), and Wasserstein-2 (light gray, dashed) distances for ({\textbf{a}}) $ Nx = 50 $ and ({\textbf{b}}) $ Nx = 100 $ triangulations.}
\label{SymDumbDists}
\end{center} 
\end{figure}

\begin{figure}[H]\begin{center}
 \begin{tabular}{cc}
 \includegraphics[width=2.5in]{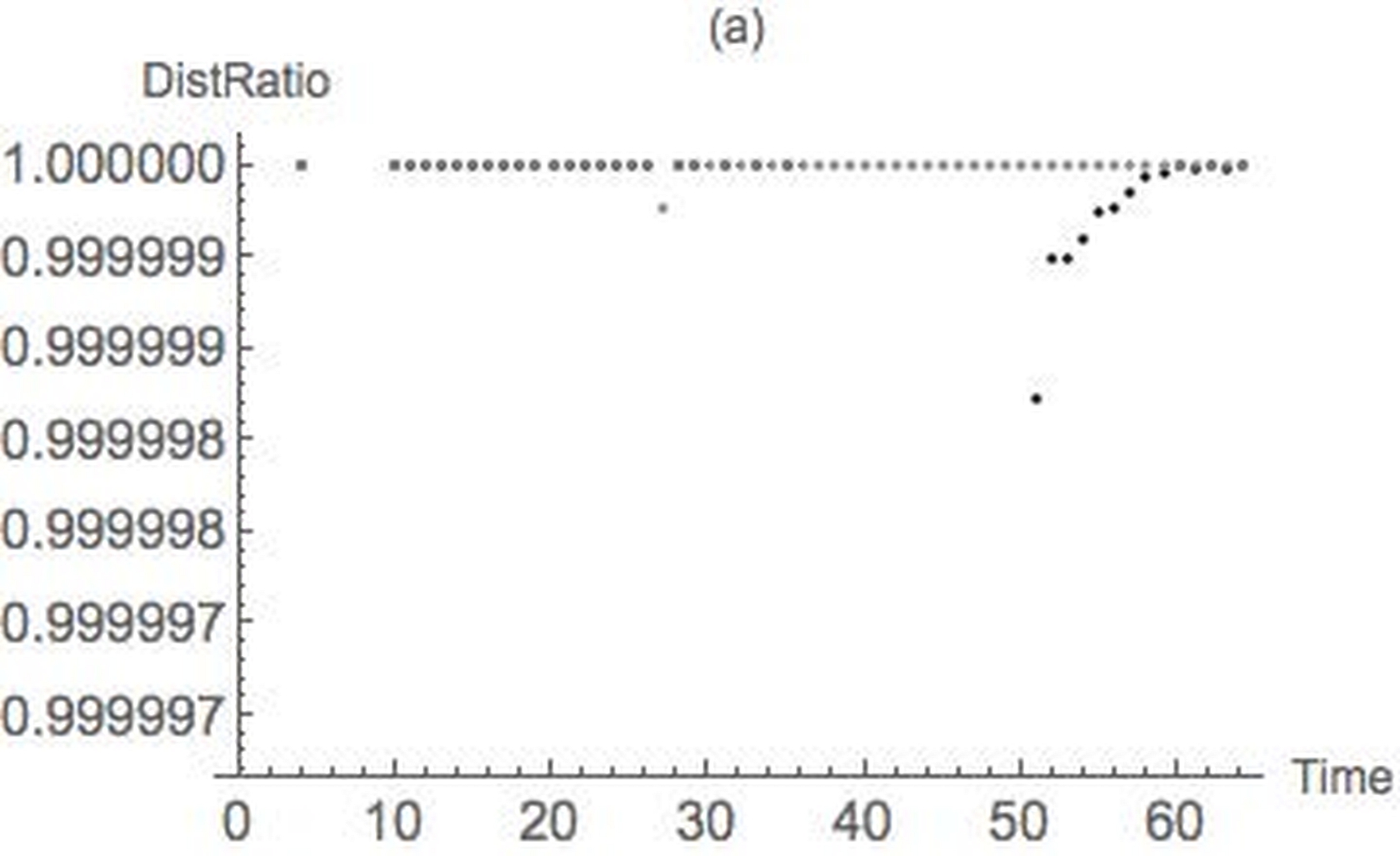} &
 \includegraphics[width=2.5in]{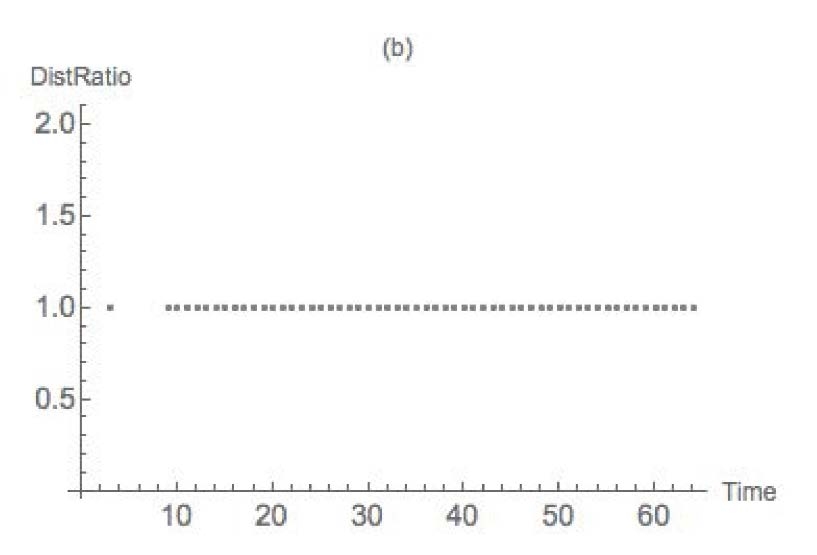}
 \end{tabular}
\caption{$ \beta_{0} $: Interpolation functions of ratios of bottleneck distance to Wasserstein-1 distance (black~dots) and Wasserstein-2 distance (gray dots) for ({\textbf{a}}) $ Nx = 50 $ and ({\textbf{b}}) $ Nx = 100 $ triangulations.}
\label{SymDumbDistRat}\end{center} 
\end{figure}

Plots of the distances (with average scalar curvature) and distance ratios for the triangulations are in Figures \ref{SymDumbDists} and \ref{SymDumbDistRat}. For most time indices of the coarse triangulation and all of them for the finer triangulation, $ d_{B} = d_{W^{1}} = d_{W^{2}} $. This occurs since, for the connected components each PD has one finite lifespan and one infinite lifespan (hence two points each, with most of the infinite lifespans being the same in adjacent PDs). For the holes, with no infinite lifespans, each PD has one point. In some instances, the difference between PDs is zero. When this occurs, the plots of the distance ratios have gaps. The differences in distance measures in the $ Nx = 50 $ triangulation are insignificant (within $ 10^{-6} $). The only differences between the respective distance measures in the $ Nx = 50 $ and $ Nx = 100 $ cases are in the orders of their magnitudes, which mirror those of the differences in lifespans between the refinements. Furthermore, the distances between adjacent PDs change little until the final four time indices. The average scalar curvature grows more steadily than the distances but has two spurts: one around time index $ 50 $ and another at the end. Regarding the distance measures, the coarser triangulation, $ Nx = 50 $, actually captures this better than the finer $ Nx = 100 $ case. This is related to whether or not certain points of the object are sampled in a given refinement. 

The distance ratios are essentially 1 at each time index (where the distance between adjacent diagrams is nonzero). From this, the samples agree that a single 
change in the geometry is occurring, with no small features contributing to changes in the geometry, and increases by several orders of magnitude from the beginning to the end of the evolution. From the distances (and~corroborated with the increase in average scalar curvature), this 
single effect occurs mainly at the end of the evolution. This can be attributed to the type of singularity formation for the symmetric dumbbell expected from RF: the singularity forms at one location on the neck.

\unskip

\subsubsection{Dimpled Dumbbell}

\begin{figure}[H]
\centering
 \includegraphics[width=5in]{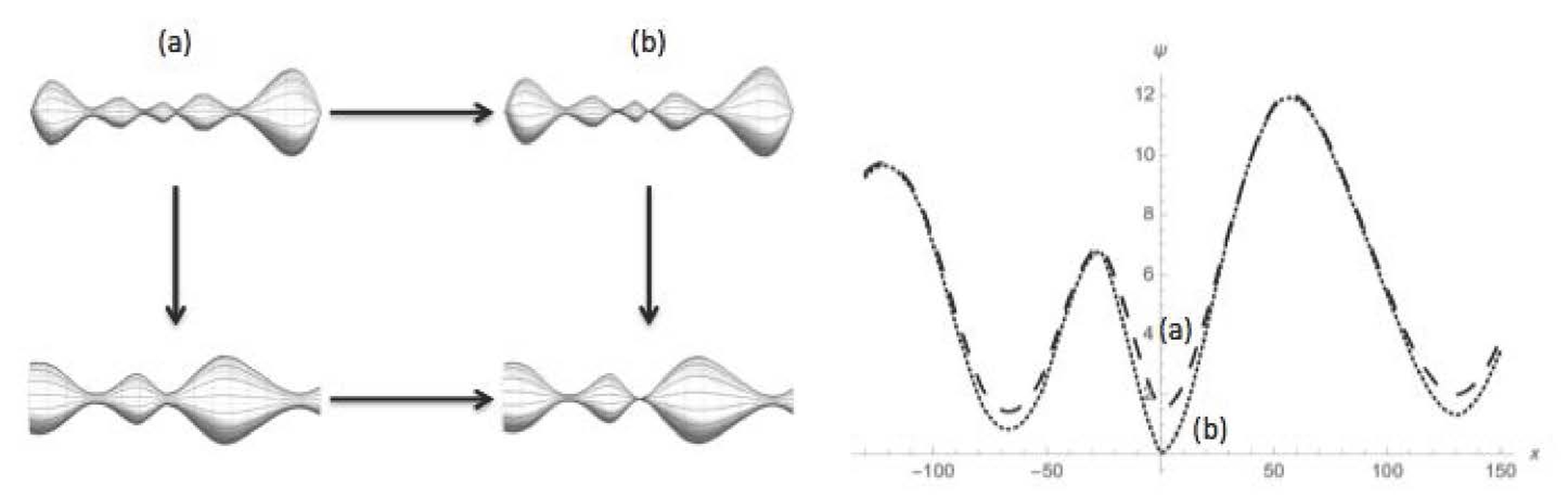}
\caption{This is an illustration of the dimpled dumbbell. The three-dimensional surfaces (above) are the plots of the dumbbell on the interval $ [-100\pi + 0.1,100\pi - 0.1] $ and close-ups of some of the necks, including the most-pinched neck, on the interval $ [-130,150] $ at initial time $ t = 0 $ ({\textbf{a}}) and final time $ t = 1.1599 $ ({\textbf{b}}). The initial (long dash) and final (short dash) radial profiles are plotted below on the interval $ [-130,150] $.}
\label{fig:DimpledDumbbell}
\end{figure}

The dimpled dumbbell is featured below (Figure \ref{fig:DimpledDumbbell}). The necks are at $ x = -180 $, $ x = -66 $, $ x = 0 $, and $ x = 130 $. The greatest pinching occurs at the neck $ x = 0 $, whose neck radius, from initial time index to final time index, 
decreases by three orders of magnitude. Furthermore, the scalar curvature at this neck grows by seven orders of magnitude. The other necks demonstrate pinching but within one order of magnitude in both neck radius and scalar curvature. Most values of $ \psi $ remain constant across the~evolution.

\begin{figure}[H]
\centering
 \begin{tabular}{cc}
  \includegraphics[width=2.5in]{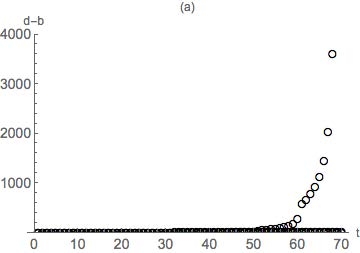} &
 \includegraphics[width=2.5in]{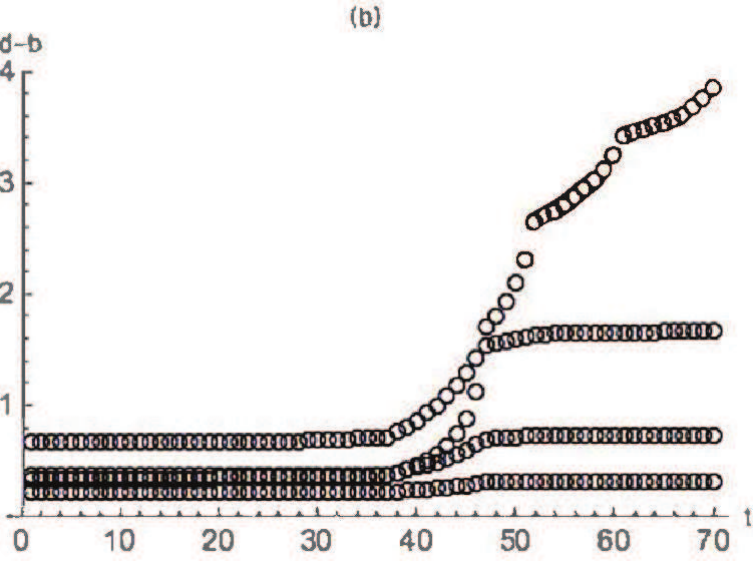}
 \end{tabular}
\caption{Lifespans $ d - b $ for $ \beta_{0} $ computed for ({\textbf{a}}) $ Nx = 25 $ and ({\textbf{b}}) $ Nx = 50 $ triangulations.}
\label{fig:GenModDumbPersistence}
\end{figure}

The triangulations are on $ Nx = 25 $ and $ Nx = 50 $ grids for 70 time indices. The lifespans for the triangulations (Figure~\ref{fig:GenModDumbPersistence}) resemble a hybrid between those of the dimpled sphere and those of the symmetric dumbbells, although this relationship is harder to see for the $ Nx = 25 $ triangulation without magnifying the shorter lifespans. In both cases, for the $ \beta_{0} $ tables, we have four finite lifespans and one infinite lifespan of low scalar curvature for all time indices. Over the evolution of the object, the $ Nx = 25 $ triangulation exhibits three finite lifespans exhibiting modest increases in birth and death values with one of the lifespans having a change in death value of eight orders of magnitude. For~the same time frame, the $ Nx = 50 $ triangulation exhibits two of the finite lifespans showing modest increases in birth and death values with the remaining lifespans featuring similar modest changes in birth values but with one order of magnitude increases in death values. 
Though $ Nx = 25 $ provides a coarser triangulation than $ Nx = 50 $, the coarser triangulation contains points in its sample with higher values of a scalar curvature at each time index. Thus, it is detecting the increase of curvature more literally than the finer triangulation. However, as is evident from Figure~\ref{fig:GenModDumbPersistence}, both triangulations capture the formation of the neckpinch singularity. Furthermore, the lifespans for the different triangulations of the dimpled dumbbell share the self-similarity of the lifespans for the different triangulations of the symmetric dumbbell.

The $ \beta_{1} $ tables have lifespans of variously short lengths but no infinite ones. Though some of the lifespans grow to lengths similar to those in the $ \beta_{0} $ tables, they are of the same order of magnitude.

\begin{figure}[H]\begin{center}
 \begin{tabular}{cc}
  \includegraphics[width=2.8in]{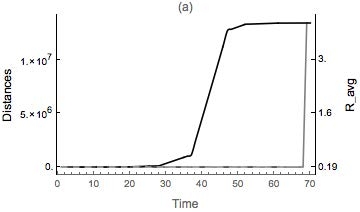} &
 \includegraphics[width=2.8in]{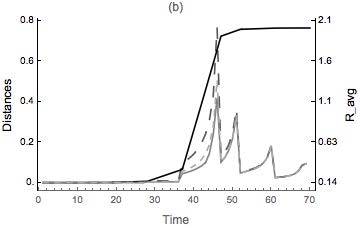}
 \end{tabular}
\caption{$ \beta_{0} $: Interpolation functions of average scalar curvature (denoted by R\_avg, solid black line), bottleneck (gray), Wasserstein-1 (dark gray, dashed), and Wasserstein-2 (light gray, dashed) distances for ({\textbf{a}}) $ Nx = 25 $ and ({\textbf{b}}) $ Nx = 50 $ triangulations.}
\label{GMPDists}\end{center}
\end{figure}

The distance measures between adjacent PDs (Figure \ref{GMPDists}) provide stronger quantitative support for this assertion. Initially, $ d_{W^{2}} $ is much closer to $ d_{B} $ than it is to $ d_{W^{1}} $ (though all are within the same order of magnitude), indicating that small changes in the geometry dominate differences in the PDs. As the final time is approached, however, $ d_{W^{1}} $ is only slightly larger than $ d_{W^{2}} $ or $ d_{B} $ (which are almost identical), indicating that a single large change dominates differences in the PDs. 

\begin{figure}[H]\begin{center}
 \begin{tabular}{cc}
 \includegraphics[width=2.5in]{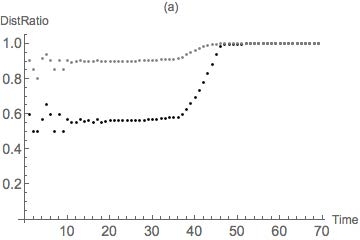} &
 \includegraphics[width=2.5in]{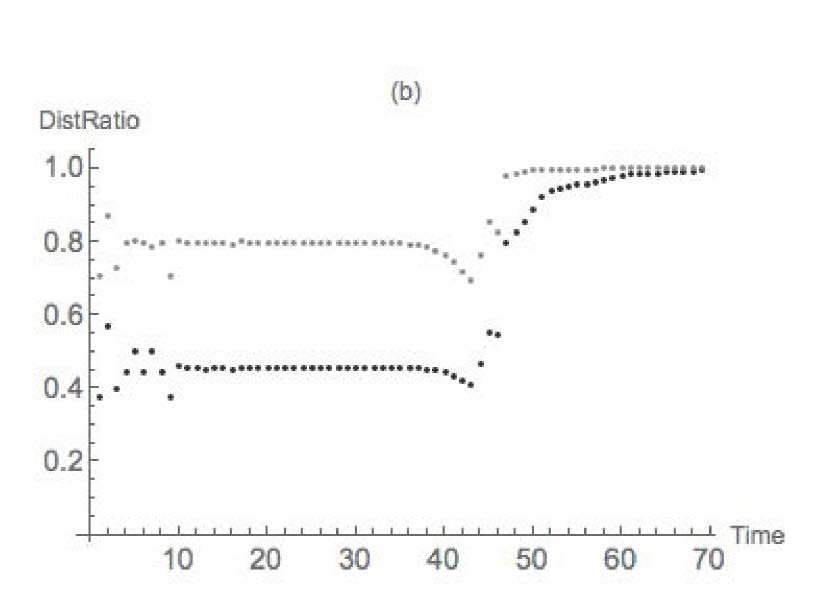}
 \end{tabular}
\caption{$ \beta_{0} $: Interpolation functions of ratios of bottleneck distance to Wasserstein-1 distance (black~dots) and Wasserstein-2 distance (gray dots) for ({\textbf{a}}) $ Nx = 25 $ and (\textbf{b}) $ Nx = 50 $ triangulations.}
\label{GMPDistRat}\end{center}
\end{figure}

The ratios (Figure \ref{GMPDistRat}) corroborate this. While $ d_{B}/d_{W^{1}} $ and $ d_{B}/d_{W^{2}} $ are initially small but oscillate, they remain respectively at around $ 0.45 $ and $ 0.8 $ until closer to the end of the evolution, when they appear to approach $ 1 $. The dimpled dumbbell shares features with the previous cases of the dimpled sphere and the symmetric dumbbell: small changes in the geometry initially dominate differences in the PDs while a single large change corresponds with a neckpinch at the end of the evolution, the~latter of which is expected from RF.

 \unskip

\section{Conclusions and Future Directions}
\label{sec:Conclusions}

The application of PH to RF provides a unique opportunity to apply computational homology to problems from differential geometry that have been theoretically developed, but in the case of the dimpled dumbbell, have not been explored numerically. Having such a ``fully-loaded'' dataset, a priori, yields a notion of how the dataset should behave. The motivation is to develop an understanding, through PH, of how a dataset that is less robust, but that describes similar systems, ought to behave. The~advantage of using PH is that it strongly connects the topological behavior of the data to the geometric behavior at each time step.

The use of two resolutions is to confirm that scale does not significantly impact the results. This is somewhat heuristic as we have examined only two resolutions for each. However, it is strongly indicative, based on the obtained numerical results, that generic effects distinguishing global and local singularities manifest regardless of the resolution.

We have available several tools for measuring the system. These are the lifespans, the distance measures between PDs and the distance ratios. The lifespans demonstrate interesting features in each of the models that, in the plots of death versus birth for the dumbbells, result in almost self-similarity. Furthermore, the values seem to be associated with the extent of clustering within the PDs at each time index. This clustering, or lack thereof, appears to be associated with the type of singularities due to the differences in the behavior of data collapsing to a round point (all points approach a common point) and neckpinch singularity formation (a single change in the data). A weakness in using only the lifespans to analyze this singularity formation is in the scale. Depending on the curvature values, the~lengths of each of the lifespans can be long, but similar.

The distances, as measured at each time index, provide a better measure. They give a notion of the weight of the contribution of small changes in the data. However, like the lifespans, they are tied to the scale of the values, and it can be difficult to assess, outside of strong contrasts in the values, the~influence of small changes. In particular, the orders of magnitudes of certain distances can obscure the notion of the closeness of values.

Of these three tools, the distance ratios provide the best measure. Because they are dimensionless, they allow for scale-independent assessments of the contribution of small features to changes in the dataset over time. The notion of closeness, then, is easier to understand than with distances alone since orders of magnitude are absorbed in the ratios.

It is interesting to note what the distance measures and distance ratios indicate between PDs near the end of the evolution for each of the models. For the dimpled sphere under RF, the singularity is global. The distances between adjacent PD's for the dimpled two-sphere do not reach a single value. Furthermore, the distance ratios fluctuate, but trend downward for both resolutions. These measurements indicate that, for the selection of time steps, many small changes are occurring between adjacent PDs. The RF evolves the dimpled two-sphere by first making the scalar curvature uniform before collapsing the object to a point, so the many small changes between the adjacent PDs coincide with the apparent merger of the points to a single one.

For the symmetric dumbbell and the dimpled dumbbell, RF produces local singularities. For~the symmetric dumbbell, the distances between adjacent PDs grow only at the end. The distance ratios are either identical or within $ 10^{-6} $ with the difference occurring in the Wasserstein-1 metric as would be expected, indicating that only a single large feature changes the geometry. For the symmetric dumbbell, only this change is occurring; it grows significantly near the end of the evolution. 

The dimpled dumbbell has a topological signal sharing features with both the dimpled sphere and the symmetric dumbbell. The distances between the PDs initially indicate that small features dominate in changing the geometry. By the end of the evolution, however, the distance measures become much closer, indicating a single large effect changing the geometry. The distance ratios confirm this as they approach one near the end of the evolution. This single change manifests mainly at the end of the evolution and coincides with the formation of a neckpinch singularity, as predicted by RF. Such~similarities of distance values for these measures indicate a single large change occurring in the data from time index to time index. 

From these results, PH is able to delineate between global and local singularities since the datasets, in these cases, will exhibit behavior that manifests in either a noticeably smaller bottleneck distance compared with the other measures or in similar values of distance measures, as verified by the distance ratios, at the end of the evolution. Similar behavior occurs at multiple resolutions indicating that PH distinguishes between global and local singularities even with a few sampled points.

For future work, we would like to investigate a variety of meshes of different sizes and shapes. The meshes considered here, while easy to implement, significantly restrict the sampling of angles and are problematic for certain curvature approximation experiments in two dimensions. Particularly, we~would like to explore tetrahedral meshes as these are ubiquitous across related disciplines, such~as~graphics.

As this work is a first step into this direction of research, we would like to consider a variety of meshes so as to enable a study of convergence. Although novel, our results are preliminary and suggest the need for further experiments involving a variety of meshes and initial data. The study of convergence requires significantly more rigor and generality than presented here. 

Regarding interdisciplinary applications, we would like to use PH to investigate systems with incomplete or sparse data. This includes problems involving the interpretation of data for hard-to-compute higher-dimensional scenarios, which may appear in the context of physical problems with more than four dimensions or systems with a large number of degrees of freedom (\mbox{e.g., phase spaces}). Furthermore, we would like to explore quantum entanglement using computational homology. The~advantage of such an approach is that one obtains sparse data from density matrices in quantum tomography; different measures of entanglement \cite{AB,CKW,WGLH} would be candidates for filtrations to fill-in the missing data as in matrix imputation. In this pursuit, it may be worthwhile to couple PH to RF in a way similar to the current work as RF is understood on continuum manifolds both real and complex, and~matrix product states (of which cluster states are) have recently been given a differential geometric representation~\cite{Osborne}. Another direction would be the extension of the methods in this work to problems when one evolves ``through'' the singularity \cite{Perelman:2003-1,Perelman:2003-2}. Furthermore, we are exploring the application of PH to the degenerate dumbbell, constructed from our interpolative dumbbell for $ \alpha \neq 1 $ and other parameters, to obtain a~homological perspective of Type-II singularity formation.

As an interesting alternative, it would be worthwhile to investigate the use of persistent homotopy~\cite{SJ} for these and other problems, especially that this approach allows flexibility in the kinds of structures considered.
\vspace{6pt} 


\section{Acknowledgments}
Howard A. Blair acknowledges support from United States Air Force (USAF) Grant \# FA8750-11-12-0275. MC acknowledges support from the University of Wisconsin-Stout Faculty Start-Up Fund and the National Academy of Sciences, National Research Council (NRC) Research Associateship program. We thank Louise Camalier, Rory Conboye and Don Sheehy for suggestions regarding this work. Warner A. Miller acknowledges support from United States Air Force (USAF) Grants \# FA8750-11-2-0089 and \# FA8750-15-2-0047 and support from the Air Force Research Laboratory/Trusted Systems Branch (AFRL/RITA) through the Summer Faculty Fellowship Program and the Visiting Faculty Research Program (VFRP) administered through the Griffiss Institute. Konstantin Mischaikow acknowledges the National Science Foundation for Grants NSF-DMS-0835621, 0915019, 1125174, 1248071 and~contracts from the Air Force Office of Scientific Research (AFOSR) and the Defense Advanced Research Projects Agency (DARPA). Vidit Nanda is grateful to Robert Ghrist for insightful discussions. His work was funded in part by the Defense Advanced Research Projects Agency Grant DARPA DSO-FA9550-12-1-0416. Any opinions, findings and~conclusions or recommendations expressed in this material are those of the author(s) and do not necessarily reflect the views of the Air Force Research Laboratory.




\end{document}